\documentclass[12pt]{amsart}%
\usepackage[T1]{fontenc}
\usepackage[utf8]{inputenc}

\usepackage{color}
\usepackage{amsfonts}
\usepackage{amsmath}
\usepackage{amssymb}
\usepackage{graphicx}
\usepackage{color}
\usepackage{enumerate}

\usepackage{amsmath,amsthm,amscd,amssymb}
\usepackage{caption}  
\usepackage{graphicx}
\usepackage{geometry} 
\numberwithin{equation}{section}

\theoremstyle{plain}
\newtheorem{theorem}{Theorem}[section]

\theoremstyle{definition}
\newtheorem{definition}[theorem]{Definition}

\newtheorem{case[theorem]}{Case}

\theoremstyle{remark}
\newtheorem{remark}[theorem]{Remark}

\numberwithin{equation}{section}

\def\hat{\widehat}

\def\be{\begin{equation}}
\def\benn{\begin{equation}\nonumber}
\def\ee{\end{equation}}
\def\bfo{\begin{eqnarray*} }
\def\efo{\end{eqnarray*} }
\def\ba{\begin{eqnarray*} }
\def\ea{\end{eqnarray*} }
\def\beq{\begin{eqnarray}}
\def\eeq{\end{eqnarray}}

\def\R{\mathbb R}

\def\hKK{\hat{K\otimes K}}
\def\hZZ{\hat{Z\times Z}}

\begin{document}

\title{Bilinear generalized Radon transforms in the plane} 


\author{A. Greenleaf, A. Iosevich, B. Krause and A. Liu}

\date{April 3, 2017}

\address{Department of Mathematics, University of Rochester, Rochester, NY, 14627}

\email{allan@math.rochester.edu}

\email{iosevich@math.rochester.edu}

\address{Department of Mathematics, University of British Columbia, Vancouver, BC, V6T 1Z2} 

\email{benkrause2323@gmail.com }

\address{Department of Mathematics, Massachusetts Institute of Technology, Cambridge, MA, 02139} 

\email{cliu568@gmail.com }

\thanks{The work of the first listed author was partially supported by NSF Grant DMS-1362271, the second by
 NSA Grant H98230-15-1-0319, and the third  by  an NSF Postdoctoral Fellowship.}

\begin{abstract} Let $\sigma$ be  arc-length measure on $S^1\subset \R^2$ and $\Theta$  denote rotation by an angle   $\theta \in (0, \pi]$. 
Define a model bilinear generalized Radon transform,
$$B_{\theta}(f,g)(x)=\int_{S^1} f(x-y)g(x-\Theta y)\, d\sigma(y),$$ 
an analogue  of the linear generalized Radon transforms  of Guillemin and Sternberg \cite{GS} and  Phong and Stein (e.g., \cite{PhSt91,St93}). 
Operators such as $B_\theta$ are motivated by problems in geometric measure theory and combinatorics. 
For $\theta<\pi$, we show that  $B_{\theta}: L^p({\Bbb R}^2) \times L^q({\Bbb R}^2) \to L^r({\Bbb R}^2)$ 
if $\left(\frac{1}{p},\frac{1}{q},\frac{1}{r}\right)\in\nolinebreak Q$,  the polyhedron 
with the vertices $(0,0,0)$, $(\frac{2}{3}, \frac{2}{3}, 1)$, $(0, \frac{2}{3}, \frac{1}{3})$, $(\frac{2}{3},0,\frac{1}{3})$, $(1,0,1)$, $(0,1,1)$ and $(\frac{1}{2},\frac{1}{2},\frac{1}{2})$, 
except for $\left( \frac{1}{2},\frac{1}{2},\frac{1}{2} \right)$,  where we obtain a restricted strong type estimate. 
For the degenerate case $\theta=\pi$, a 
more restrictive set of exponents holds.
In the scale of normed spaces, $p,q,r \ge 1$, the type set $Q$ is sharp. 
Estimates for the same exponents are also proved for a class of bilinear generalized Radon transforms in $\R^2$ of the form 
$$ B(f,g)(x)=\int \int \delta(\phi_1(x,y)-t_1)\delta(\phi_2(x,z)-t_2) \delta(\phi_3(y,z)-t_3) f(y)g(z) \psi(y,z) \, dy\, dz, $$ 
where $\delta$ denotes the Dirac distribution, $t_1,t_2,t_3\in\R$, $\psi$ is a smooth cut-off  and the defining functions $\phi_j$ satisfy some natural geometric assumptions. 
\end{abstract} 

\maketitle

\section{Introduction}

\vskip.125in 

A classical result due to Littman \cite{L71} and Strichartz \cite{Str70} (see also \cite{O04}) says that, for $d\ge 2$,  
the spherical averaging  (or F. John  \cite{J}) operator,
\begin{equation} \label{sa} 
Af(x)=\int_{S^{d-1}} f(x-y)\, d\sigma(y)=\sigma*f(x), 
\end{equation} 
where $\sigma$ is the surface measure on $S^{d-1}$, is bounded from $L^p({\Bbb R}^d)$ to $L^q({\Bbb R}^d)$ 
iff $\left(\frac{1}{p}, \frac{1}{q} \right)$ is in the closed triangle with the vertices 
\begin{equation} \label{linearrange} 
(0,0), (1,1), \ \text{and} \ \left(\frac{d}{d+1}, \frac{1}{d+1}\right). 
\end{equation}

The operator $A$ is, along with the classical Radon transform,  a model for the generalized Radon transforms studied by Guillemin and Sternberg \cite{GS} and Phong and Stein  (see, e.g., \cite{PhSt91,St93}, and the references contained therein). These are  linear operators of the form 
\begin{equation} \label{radontransformgen} {\mathcal R}f(x)=\int_{\phi(x,y)=t} f(y) \psi(x,y) \, d\sigma_{x,t}(y), \end{equation} 
where $t\in\R$, $\phi\in C^\infty(\R^d\times\R^d)$ is a defining function, i.e., $d_{x,y}\phi\ne (0,0)$ on 
$$Z_t:=\left\{(x,y)\in\R^d\times\R^d:\phi(x,y)=t\right\},$$
$\psi$ is a smooth cut-off, 
$\sigma_{x,t}$ is  surface measure on $Z_t$,  
and  $\phi$ satisfies   the following condition \cite{PhSt91}:

\begin{definition} 
A defining function $\phi: {\Bbb R}^d \times {\Bbb R}^d \to {\Bbb R}$ satisfies the \emph{Phong-Stein rotational curvature condition} at $t$ if,
for all  $(x,y)\in Z_t$,
\begin{equation} \label{cond PS}
det
\begin{pmatrix} 
 0 & \nabla_{x}\phi \\
 -{(\nabla_{y}\phi)}^{T} & \frac{\partial^2 \phi}{\partial x_i \partial y_j}
\end{pmatrix}
\neq 0.
\end{equation} 
\end{definition} 

Under the rotational curvature assumption, ${\mathcal R}: L^p({\Bbb R}^d) \to L^q({\Bbb R}^d)$ for $\left(\frac{1}{p}, \frac{1}{q}\right)$ as in (\ref{linearrange}) above. 
This is a folk theorem (as far as we know), and follows by substituting the $L^2\to L^2_{\frac{d-1}2}$ boundedness of Fourier integral operators associated with canonical graphs  into Strichartz's proof \cite{Str70} in the case of the spherical averaging operator. 
Note that if $\phi(x,y)=|x-y|$, the Euclidean distance, and $t=1$, we recover the spherical averaging operator $A$ of (\ref{sa}). 
\medskip

The purpose of this paper is to study natural bilinear variants of  the linear generalized Radon transforms, with the considerations   limited to two dimensions. 
A family of model operators, arising from combinatorial geometry,  is given by 
\begin{equation} \label{bisa} 
B_{\theta}(f,g)(x)=\int f(x-y)g(x-\Theta y) \, d\sigma(y), 
\end{equation} where $\sigma$ 
is the arc-length measure on $S^1$ and 
$\Theta$ denotes the counter-clockwise rotation by an angle $\theta\ne 0$. 
We exclude the degenerate case $\theta=0$, since  
$B_0(f,g)=A(f\cdot g)$, with the linear circular mean operator $A$ as in (\ref{sa}).

Before stating the main theorem,
we describe some motivating applications. Consider $n$ points in a point set $P\subset\R^2 $ and  the problem of counting equilateral triangles of side-length $1$ among  them (see Fig. 1 above for a particularly triangle-rich configuration). We have 
$$ \# \{(x,y,z) \in P^3: |x-y|=|x-z|=|y-z|=1 \}$$
\begin{equation} \label{triform}\qquad\quad =\sum_{x,y,z} 1_C(x-y) 1_C(x-z)1_C(y-z) 1_P(x)1_P(y)1_P(z), \end{equation} 
where $C$ is the circle of radius $1$ centered at the origin. This expression equals 
$$ \sum_{x,u,v} K(u,v) 1_P(x-u) 1_P(x-v) 1_P(x),$$ where $K$ is the indicator function of the set 
$$ \{(u,v) \in C \times C: |u-v|=1\}=\{(u, \Theta u): u \in C \} \cup \{(u, \Theta^{-1}u): u \in C\},$$ 
where $\Theta$ is the rotation by $\frac{\pi}{3}$. 

\begin{figure}
\label{trianglelattice}
\centering
\includegraphics[scale=.5]{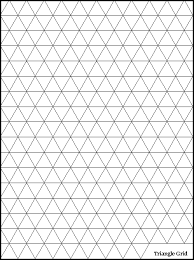}
\caption{An equilateral  triangular grid (mathforum.org)}
\end{figure}

It follows that the trilinear form in (\ref{triform}) equals 
$$ \sum_x {\mathcal B}_{\frac\pi3}(1_P,1_P)(x) 1_P(x)+\sum_x {\mathcal B}_{-\frac\pi3}(1_P(x),1_P(x)) 1_P(x),$$ 
where ${\mathcal B}_{\theta}$ is the discrete version of the bilinear operator defined in (\ref{bisa}) above. Different values of $\theta$ are similarly associated with counting triangles of different congruence types. See, for example, \cite{B86,FKW90,GI12} where operators of this types are studied, in one form or another, in the context of point configuration problems in geometric measure theory. 

The spherical averaging operator can be similarly interpreted as the continuous analogue of an operator counting pairs of distances. Indeed, arguing as above, let $P$ be a finite point set, let $C$ be the unit circle and consider 
$$ \# \{(x,y) \in P \times P: |x-y|=1 \}=\sum_{x,y} 1_C(x-y) 1_P(x)1_P(y)$$
$$=\sum_x \left\{ \sum_y 1_P(x-y) 1_C(y) \right\} 1_P(x)=\sum_x ({\mathcal A}1_E)(x) 1_E(x),$$ where ${\mathcal A}$ is  the discrete analogue of the spherical averaging operator $Af(x)$. 

A way to understand the operators $Af(x)$ and $B_{\theta}(f,g)(x)$ in terms of a coherent geometric paradigm is the following. Let $E$ be a compact subset of ${\Bbb R}^d$. Define a graph 
by designating the vertices to be  the points of $E$,  and connecting two vertices $x$ and $y$ by an edge iff $|x-y|=1$. Then the spherical averaging operator $Af(x)$ may be viewed as the 
edge operator on this graph. Now define a hyper-graph on $E$ by connecting a triple $x,y,z$ by a hyper-edge iff $|x-y|=|x-z|=|y-z|=1$. 
The hyper-edge operator  is then precisely the bilinear 
operator  $B_{\frac\pi3}(f,g)$ (or $B_{-\frac{\pi}{3}}(f,g)$). 
One can define similar objects by replacing the distance function $|x-y|$ with a more general function $\phi(x,y)$, as in (\ref{radontransformgen}). 
\medskip

These examples suggest a natural class of bilinear Radon transforms  of the form 
\begin{equation} \label{generaloperator} 
B(f,g)(x)=\lim_{\epsilon \downarrow 0} \epsilon^{-3} \int \int_{|\phi_1(x,y)-t_1|<\epsilon; |\phi_2(x,z)-t_2|<\epsilon; |\phi_3(y,z)-t_3|<\epsilon} \psi(y,z)f(y)g(z) dydz, \end{equation} 
where $\psi$ is a smooth cut-off function and $\phi_j$'s are suitably regular functions. In the case when $\phi_j(x,y) \equiv |x-y|$, we recover the operator $B_{\theta}$, $\theta=\frac{\pi}{3}$ 
defined in (\ref{bisa}). While there has been considerable progress in the study of bilinear analogues of singular integral  and pseudodifferential operators,
e.g., by Coifman-Meyer \cite{CM90}, Demeter-Tao-Thiele \cite{DTT08}, Grafakos-Li \cite{GL04}, Lacey-Thiele \cite{LT97}, Muscalu-Tao-Thiele \cite{MTT02}, 
Grafakos-Torres \cite{GTr02} and others,   bilinear generalized Radon transforms has not yet been widely studied. 
We take a small step in this direction in the current paper. 

\vskip.125in 

Our result for the model operators $B_\theta$ is the following; an extension to a more general class of bilinear generalized Radon transforms is in Section \ref{generaltransforms} below. 

\vskip.125in 

\begin{theorem} \label{main} 
For $\theta\in (0,2\pi)$, let $B_{\theta}$ be defined as in (\ref{bisa}) above. 

(i) Suppose that $\theta \not= \pi$. Then the  type set $Q(B_{\theta})$ of $B_{\theta}$, i.e., those $\left( \frac{1}{p}, \frac{1}{q}, \frac{1}{r} \right)\in [0,1]^3$ such that
$$B_{\theta}: L^p({\Bbb R}^2) \times L^q({\Bbb R}^2) \to L^r({\Bbb R}^2),$$ 
is the closed polyhedron  with vertices $(0,0,0)$, $(\frac{2}{3}, \frac{2}{3}, 1)$, $(0, \frac{2}{3}, \frac{1}{3})$, $(\frac{2}{3},0,\frac{1}{3})$, $(1,0,1)$, $(0,1,1)$ 
and $(\frac{1}{2},\frac{1}{2},\frac{1}{2})$,  
except for $(\frac{1}{2}, \frac{1}{2}, \frac{1}{2})$, where  a restricted strong type bound holds, i.e., $B_\theta:L^{2,1}\times L^{2,1}\to L^2$. (See Fig. 2 below.)
\medskip

(ii) If $\theta=\pi$,  the operator is bounded if $\left( \frac{1}{p}, \frac{1}{q}, \frac{1}{r} \right)\in Q(B_{\pi})$,  the closed polyhedron with vertices  $(0,0,0)$, 
$(\frac{2}{3}, \frac{2}{3}, 1)$, $(0, \frac{2}{3}, \frac{1}{3})$, $(\frac{2}{3},0,\frac{1}{3})$, $(0,1,1)$, and $(1,0,1)$. (See Fig. 3 below.) 
\medskip

(iii) Moreover, in  the Banach cube $p,q,r\ge 1$, the exponents in both (i) and (ii)  are best possible, 
except for the question of whether in case (i)  there is a strong type estimate for $(\frac{1}{2}, \frac{1}{2}, \frac{1}{2})$, which is unknown at this time. 
\end{theorem} 

\begin{remark} By reflection about the horizontal axis, it suffices to consider $0<\theta\le \pi$.
\end{remark}

\vskip.125in 

\begin{remark} Examination of the proof shows that all the non-trivial endpoints in the degenerate case, $\theta=\pi$, follow, 
in one way or another, from the $L^p({\Bbb R}^2) \to L^q({\Bbb R}^2)$ bounds (\ref{linearrange}) for the circular averaging operator $f\to Af$ defined in (\ref{sa}) above. 
In the non-degenerate case, there is restricted strong type boundedness at an additional vertex, $\left(\frac{1}{2}, \frac{1}{2}, \frac{1}{2}\right)$, 
and the estimates which follow from that by interpolating with the estimates valid in the degenerate case.. 
\end{remark} 

\vskip.125in 

\begin{remark} If we do not restrict ourselves to the Banach cube $p,q,r\ge 1$, the sharpness examples suggest the possibility of an additional estimate,
$$ B_{\theta}: L^1({\Bbb R}^2) \times L^1({\Bbb R}^2) \to L^{\frac{1}{2}}({\Bbb R}^2),$$ 
but we do not know  whether or not this bound holds. 
\end{remark} 

\vskip.125in

\vskip.25in 

\section{A general class of  bilinear generalized Radon transforms} 
\label{generaltransforms} 

\vskip.125in

We now extend Theorem \ref{main} to a class of bilinear operators modeled on the linear generalized Radon transforms of \cite{GS,PhSt91}.
Bilinear generalized Radon transforms (and more singular variants) have arisen in inverse problems (e.g., \cite{CFGGN,GLSSU}), 
although in those problems the relevant estimates are in terms of Sobolev spaces.
\vskip.125in

 Let $(\phi_1,\phi_2,\phi_3)$ be a trio of smooth, real-valued  functions on $\R^2\times\R^2$,
and for $(t_1,t_2,t_3)\in \R^3$, 
let $K$ be the product delta distribution on $\R^2\times\R^2\times\R^2$,
\be\label{def K}
K(x,y,z)=\delta(\phi_1(x,y)-t_1)\cdot\delta(\phi_2(x,z)-t_2)\cdot\delta(\phi_3(y,z)-t_3).
\ee
We define a bilinear generalized Radon transform by
\be\label{def B}
B(f,g)(x)=\int\int K(x,y,z)\, f(y)g(z)dydz,
\ee
defined weakly by
\be\label{def B weak}
\langle B(f,g),h\rangle=\langle K,h(x)f(y)g(z) \rangle.
\ee

\noindent If one defines 
$$Z_1=\{ (x,y)\in \R^2\times\R^2: \phi_1(x,y)=t_1\},\quad Z_2=\{(x,z)\in \R^2\times\R^2: \phi_2(x,z)=t_2\},$$

$$\hbox{ and } Z_3=\{(y,z)\in \R^2\times\R^2: \phi_3(y,z)=t_3\},$$
then $Z_1,\, Z_2$ and $Z_3$ are smooth submanifolds in $\R^2\times\R^2$ under the assumption that
the $\phi_j$ are \emph{defining functions} in the sense that
\be\label{lin indep}
d\phi_j\ne (0,0)\hbox{ on }Z_j,\, j=1,2,3.
\ee
We now slightly strengthen this assumption: Setting
$$Z=\{(x,y,z)\in \R^2\times\R^2\times\R^2: (x,y)\in Z_1, (x,z)\in Z_2, (y,z)\in Z_3\},$$
which is the support of $K(x,y,z)$, it will follow that $Z\subset \R^6$ is a smooth, codimension 3 submanifold,
and the product in (\ref{def K}) is well-defined,
if $\left(\phi_1(x,y),\phi_2(x,z),\phi_3(y,z)\right)$ is a set of defining functions for $Z$, i.e., one has 
\be\label{smooth Z}
 rank\left[
\begin{matrix}
d_x\phi_1 & d_y\phi_1& 0 \cr
d_x\phi_2 & 0 & d_z\phi_2\cr
0 & d_y\phi_3 & d_z \phi_3
\end{matrix}
\right]=3,\quad \forall (x,y,z)\in Z.
\ee
$K$ is then a conormal distribution with respect to $Z$ of order 0, denoted $K\in I^0(Z)$, since it has a representation
\be\label{osc rep K}
K(x,y,z)=\int_{\R^3}e^{i[\phi_1(x,y)\tau_1+\phi_2(x,z)\tau_2+\phi_3(y,z)\tau_3]} 1(\tau)\, d\tau,
\ee 
interpreted as an oscillatory integral in the sense of H\"{o}rmander \cite{Hor}.
Equivalently, $K$ is a Fourier integral distribution associated with the Lagrangian manifold $N^*Z\subset T^*\R^6\setminus 0$, 
$K\in I^0(N^*Z)$, 
where $N^*Z$ denotes the conormal bundle of $Z$,
\beq\label{conormal bundle}
N^*Z=& &\Big\{(x,y,z,\xi,\eta,\zeta): (x,y,z)\in Z,\, \xi=\tau_1d_x\phi_1+\tau_2 d_x\phi_2,\, \eta=\tau_1 d_y\phi_1+\tau_3 d_y\phi_3,\nonumber\\
& &\qquad\qquad\qquad\qquad  \zeta=\tau_2d_z\phi_2+\tau_3d_z\phi_3,\quad (\tau_1,\tau_2,\tau_3)\in\R^3\setminus 0\Big\}.
\eeq
The order of $K$ as a conormal distribution just equals the order, 0, of its amplitude $1(\tau)$, while its order as a Fourier integral distribution  is determined by $0+\frac32-\frac64=0$.
\bigskip

The following is the extension of Theorem \ref{main} to this more general class of bilinear generalized Radon transforms.
The second order rotational curvature condition (\ref{cond PS}) was discussed in the Introduction, 
while the first order conditions  (\ref{cond yz'}) - (\ref{cond yy'})  on the defining functions and the surface $\hZZ$  will be defined in Sec. \ref{sec L2 general} below.
\medskip

\begin{theorem}\label{thm general}
A  bilinear generalized Radon transform $B$ defined by (\ref{def B}) with kernel (\ref{def K}) conormal for a smooth $Z\subset\R^2\times\R^2\times\R^2$ with defining functions $\phi_1,\phi_2,\phi_3$ satisfying (\ref{smooth Z}) has the same type set  as the non-degenerate model operators $B_\theta$ in Theorem \ref{main} (i) if 
\medskip

(a) $\phi_3$ satisfies the  rotational curvature condition (\ref{cond PS})
\smallskip

\noindent and
\smallskip

(b) $\forall (x,y,z,y',z')\in \hZZ$,  at least one of (\ref{cond yz'}) or (\ref{cond zy'}) or [ (\ref{cond zz'}) and (\ref{cond yy'}) ] holds.

\end{theorem}


\section{Sharpness examples} 

\vskip.125in 

In this section we show that the ranges of exponents  in Theorem \ref{main} are best possible. 

\subsection{General rotation $\theta$:} Let $f(x)=\chi_{B_{\delta}}(x)$ and $g(x)=\chi_{A_{\delta}}(x)$, where $A_{\delta}$ is the annulus of radius $2$ and width $\delta$. Then $B_{\theta}(f,g)(x) \approx \delta$ on the annulus of radius $1$ and width $\delta$. It follows that 
$$ {||B_{\theta}(f,g)||}_r \approx \delta^{1+\frac{1}{r}},$$ which leads to the relation 
$$ \delta^{1+\frac{1}{r}} \lesssim \delta^{\frac{2}{p}} \cdot \delta^{\frac{1}{q}}, $$ or, equivalently 
\begin{equation} \label{ballanusleft} \frac{2}{p}+\frac{1}{q} \leq 1+\frac{1}{r}. \end{equation} 

By symmetry we also obtain 
\begin{equation} \label{ballanusright} \frac{2}{q}+\frac{1}{p} \leq 1+\frac{1}{r}. \end{equation} 

Taking $p=q$ yields 
$$ \frac{3}{p} \leq 1+\frac{1}{r}.$$

Taking $r=1$, this yields the condition $p \ge \frac{3}{2}$. When $r=2$, we get $p \ge 2$. 

\vskip.125in 

Let $E=F=B_R$, the ball of radius $R$ centered at the origin. It follows that $B_{\theta}(\chi_E,\chi_F)(x) \approx 1$ on a set of measure $\approx R^2$. It follows that 
$$ R^{\frac{2}{r}} \lesssim R^{\frac{2}{p}} \cdot R^{\frac{2}{q}},$$ which leads to the restriction 
\begin{equation} \label{largeballexample} \frac{1}{r} \leq \frac{1}{p}+\frac{1}{q}. \end{equation} 

\vskip.125in 

\subsection{The case $\theta=\pi$:} 

Let $E$ be the indicator function of the $\epsilon$ by $\epsilon^2$ rectangle tangent to $S^1$ at the north pole and let $F$ be the same object at the south pole. Then $B(\chi_E, \chi_F)(x) \approx \epsilon$ on a set of area $\approx \epsilon^3$. It follows that 
$$ \epsilon \cdot \epsilon^{\frac{3}{r}} \leq C \epsilon^{\frac{3}{p}+\frac{3}{q}},$$ which forces 
\begin{equation} \label{2rectangles} \frac{3}{p}+\frac{3}{q} \leq 1+\frac{3}{r}, \end{equation} which is a stricter condition than (\ref{ballanusleft}) or (\ref{ballanusright}). 

\vskip.125in 

\subsection{The case $\theta \not=\pi$} 

Let $E$ be the indicator function of the $\epsilon$ by $\epsilon^2$ rectangle tangent to $S^1$ at the north pole and let $F$ be the same object tangent to $S^1$ at the angle $\theta \not=0, \pi$. Then $B(\chi_E, \chi_F)(x) \approx \epsilon$ on a set of area $\approx \epsilon^4$. It follows that 
$$ \epsilon \cdot \epsilon^{\frac{4}{r}} \leq C \epsilon^{\frac{3}{p}+\frac{3}{q}},$$ which forces 
\begin{equation} \label{2rectanglesnondeg} \frac{3}{p}+\frac{3}{q} \leq 1+\frac{4}{r}, \end{equation} which is, once again a  stricter condition than (\ref{ballanusleft}) and (\ref{ballanusright}). 

\vskip.125in 

\subsection{Sharpness examples  from duality} 

\vskip.125in 

Fixing a function $g$, define $T_1f(x)=B_{\theta}(f,g)(x)$. Then 
$$ <T_1f,h>=\int \int f(x-y)g(x-\Theta y)h(x)d\sigma(y)dx.$$ 
Let $x'=x-y$, one obtains
$$ \int \int f(x') g(x'+y-\theta y)h(x'+y) dx' d\sigma(y)=<h,T_1^{*}f>,$$ where 
$$ T_1^{*}h(x')=\int h(x'+y) g(x'+y-\Theta y) d\sigma(y)=:B^1_{\theta}(h,g)(x').$$ 
Similarly, with $f$  fixed, let $T_2g(x)=B_{\theta}(f,g)(x)$. Then 
$$ T_2^{*}h(x')=\int h(x'+\Theta y) f(x'+\Theta y-y) d\sigma(y)=:B^2_\theta(f,h)(x').$$ 

It is not difficult to see that $T_j^{*}$, $j=1,2$, satisfies the same bounds $T$ does. 
This is because $T_j^{*}$ has essentially the same form with respect to another curve with strictly positive curvature. 
In other words, if $Q(B_{\theta})$ denotes the type set of triples $\left(\frac{1}{p}, \frac{1}{q}, \frac{1}{r} \right)$ 
such that $B_{\theta}: L^p({\Bbb R}^2) \times L^q({\Bbb R}^2) \to L^r({\Bbb R}^2)$,
then, for both $j=1,2$,  $\left(\frac{1}{p}, \frac{1}{q}, \frac{1}{r} \right) \in Q(B_{\theta})$ if and only if $\left(\frac{1}{p}, \frac{1}{q}, \frac{1}{r} \right) \in Q(B^j_{\theta})$. 
\medskip

Applying this idea to (\ref{ballanusleft}) and (\ref{ballanusright}), we obtain the constraint 
\begin{equation} \label{ballanusdual} \frac{1}{p}+\frac{1}{q} \leq \frac{2}{r}. \end{equation} 

On the other hand, applying duality to (\ref{2rectanglesnondeg}) yields constraints
\begin{equation} \label{rectdualleft} \frac{4}{p}+\frac{3}{q} \leq 2+\frac{3}{r} \end{equation} and 
\begin{equation} \label{rectdualright}  \frac{3}{p}+\frac{4}{q} \leq 2+\frac{3}{r}. \end{equation} 

\vskip.125in 

\section{Summary of sharpness conditions and vertices of $Q(B_{\theta})$} 

\vskip.125in 

The following is the list of exponent restrictions in the setting of Banach spaces.  
\bigskip

\begin{itemize} 

\item $ 0 \leq \frac{1}{p}, \frac{1}{q}, \frac{1}{r} \leq 1, \ (\text{Banach cube}). $

\item $ \frac{2}{p}+\frac{1}{q} \leq 1+\frac{1}{r},  \ \frac{1}{p}+\frac{2}{q} \leq 1+\frac{1}{r}, \ (\text{Small ball and annulus}). $

\item $ \frac{1}{p}+\frac{1}{q} \leq \frac{2}{r}, (\text{Dual of small ball and annulus}).$

\item $ \frac{3}{p}+\frac{3}{q} \leq 1+\frac{4}{r}, (\text{Tangent rectangles} \ \theta\not=\pi).$

\item $ \frac{3}{p}+\frac{3}{q} \leq 1+\frac{3}{r}, (\text{Tangent rectangles} \ \theta=\pi).$

\item $ \frac{4}{p}+\frac{3}{q} \leq 2+\frac{3}{r}, \  \frac{3}{p}+\frac{4}{q} \leq 2+\frac{3}{r}, \ (\text{Dual of tangent rectangles} \ \theta \not=0,\pi).$ 

\item $  \frac{1}{r} \leq \frac{1}{p}+\frac{1}{q}, \ \text{(Large ball)}.$ 

\end{itemize} 

\vskip.125in 

\begin{remark} Boxes with dimension $\delta \times \delta \times \dots \times \delta \times \delta^2$, tangent to the unit sphere, are often referred to as C. Fefferman boxes \cite{F73} in the harmonic analysis literature. \end{remark} 

\vskip.125in 

\subsection{Vertices of $Q(B_{\theta})$} \label{verticeslist} 
Using SageMath \cite{Sage} one can compute the vertices of the polyhedron determined by the inequalities above. 
In the case $\theta \not=\pi$, the vertices are
\bigskip 

 (i) $(0,0,0)$
 
 (ii) $(0, 1, 1)$
 
 (iii) $(1, 0, 1)$ 
 
 (iv) $\left(\frac{2}{3}, 0, \frac{1}{3}\right)$
 
 (v) $\left(0, \frac{2}{3}, \frac{1}{3}\right)$
 
 (vi) Universal vertex: $\left(\frac{2}{3}, \frac{2}{3}, 1\right)$

 (vii) Non-degenerate vertex: $\left(\frac{1}{2}, \frac{1}{2}, \frac{1}{2} \right)$

\vskip.125in

In the case $\theta=\pi$, the vertices are 
\bigskip

 (i) $(0,0,0)$ 
 
(ii) $(0, 1, 1)$

(iii) $ (1, 0, 1)$

 (iv) $ \left(\frac{2}{3}, 0, \frac{1}{3}\right)$
 
(v) $\left(0, \frac{2}{3}, \frac{1}{3}\right)$
 
(vi) Universal vertex: $\left(\frac{2}{3}, \frac{2}{3}, 1 \right)$

\vskip.125in 

\begin{remark} We refer to $\left(\frac{2}{3}, \frac{2}{3}, 1 \right)$ as the ``universal vertex" because it arises for every $\theta$. 
We refer to $\left(\frac{1}{2}, \frac{1}{2}, \frac{1}{2} \right)$ as the non-degenerate vertex because it only arises in the case $\theta \not=\pi$. 
Note that the vertices (i)-(vi) are the same for both the non-degenerate and  degenerate cases. 

\end{remark} 

\vskip.125in 

\begin{figure}
\label{nondegenerate}
\centering
\includegraphics[scale=.6]{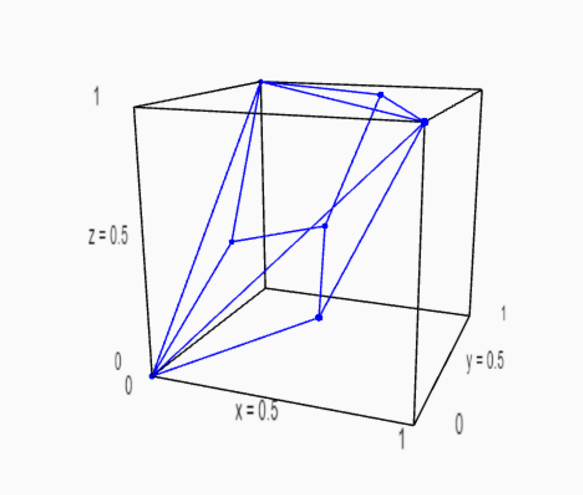}
\caption{A plot of the typeset  $ Q(B_{\theta})$ for the non-degenerate cases,  $\theta \not=\pi$.}
\end{figure}

\begin{figure}
\label{degenerate}
\centering
\includegraphics[scale=.6]{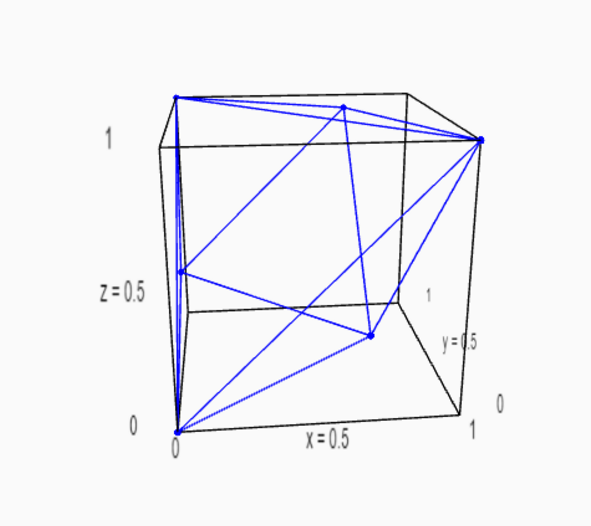}
\caption{A plot of the typeset  $Q(B_{\pi})$ in the degenerate case, $\theta=\pi$. }
\end{figure}

\vskip.125in


\section{Trivial bounds} 

We are going to establish boundedness of $B_\theta$ at the vertices described above. The full range of exponents is then recovered using multi-linear interpolation; see, e.g., \cite{GT03} and the references contained therein. 
One may assume throughout that $f,g \ge 0$ , since the general case  can be recovered  by writing $f$ and $g$ 
in terms of their real and imaginary parts, and then these as differences of their positive and negative parts.
We have the pointwise estimate
$$ |B_{\theta}(f,g)(x)|= \int f(x-y) g(x-\Theta y) d\sigma(y) \leq ||f||_{L^\infty}\cdot ||g||_{L^\infty},$$ 
hence $(0,0,0)$, the vertex (i) in Subsec. \ref{verticeslist}, is in the simplex of exponents where $B_{\theta}$ is bounded. 
\medskip

Similarly, 
$$ |B_{\theta}(f,g)(x)| \leq \left(\int f(x-y) d\sigma(y)\right)\cdot ||g||_{L^\infty}$$ and 
$$ |B_{\theta}(f,g)(x)| \leq ||f||_{L^\infty}\cdot\left(\int g(x-\Theta y) d\sigma(y)\right),$$ 
so that
$$ B_{\theta}: L^{\infty}({\Bbb R}^d) \times L^q({\Bbb R}^2) \to L^r({\Bbb R}^d) \ \text{and} \ B_{\theta}: L^p({\Bbb R}^2) \times L^{\infty}({\Bbb R}^2) \to L^r({\Bbb R}^2)$$ for $(\frac{1}{q}, \frac{1}{r})$ and 
$(\frac{1}{p}, \frac{1}{r})$ in the  typeset for the circular averaging operator, $A$,
which is the triangle with vertices $(0,0)$, $(1,1)$ and $(\frac{2}{3}, \frac{1}{3})$ (i.e., (\ref{linearrange}) for $d=2$). 
This proves boundedness of $B_\theta$ at the vertices (ii)-(v) of Subsec. \ref{verticeslist} for all $\theta\in (0,\pi]$. 

\vskip.125in

\section{The $L^{\frac{3}{2}} \times L^{\frac{3}{2}} \to L^1$ estimate for the model operators} 
\label{strichartz}

\vskip.125in 

For $f,g\ge 0$, writing $\int B_\theta(f,g)(x)\, dx$ as
$$ \int \int f(x-y) g(x-\Theta y) d\sigma(y) dx
=\int f(x) \left\{ \int g(x+y-\Theta y) d\sigma(y) \right \} dx
=:\int f(x) \cdot (A_{\theta}g)(x) dx,$$ 
and applying Holder, we see that 
\begin{equation} \label{vzhopu} 
||B_\theta(f,g)||_{L^1(\R^2)}\leq {||f||}_{L^{\frac{3}{2}}({\Bbb R}^2)} \cdot {||A_{\theta}g(x)||}_{L^{3}({\Bbb R}^2)}. 
\end{equation}
Observe that 
$$ {|y- \Theta y|}^2=2(1-<y, \Theta y>),$$ which for $|y|=1$ is a  constant depending only on $\theta$, non-zero provided that 
$\Theta\ne I$, the identity map. 
Hence, $A_\theta$ is a rescaled version of the circular averaging operator $A$ from (\ref{sa}) and satisfies the same estimates
(\ref{linearrange}), in particular the $L^{\frac32}\to L^3$ bound.
 Therefore, if $\Theta \not=I$, 
$$ ||B_\theta(f,g)||_{L^1(\R^2)}\lesssim {||f||}_{L^{\frac{3}{2}}({\Bbb R}^2)} {||g||}_{L^{\frac{3}{2}}({\Bbb R}^2)}$$ 
by the classical result of Strichartz \cite{Str70} and Littman \cite{L71}. 
This establishes boundedness of $B_\theta$ at the vertex (vi) in  Subsection \ref{verticeslist}  for all $\theta\in (0,\pi]$. 

\vskip.25in

\section{The $L^{2,1} \times L^{2,1} \to L^2$ estimate for the model operator, $\theta \not=\pi$} 
\label{sec L2 model}

\vskip.125in 

We want to show that $B_\theta$ is of restricted strong type,  that is, $B_\theta:L^{2,1}({\Bbb R}^2) \times L^{2,1}({\Bbb R}^2) \to L^2({\Bbb R}^2)$,
for $\theta\ne\pi$. Thus, we need to show that if $E,\, F\subset\R^2$, then  $||B_\theta(\chi_E,\chi_F)||_{L^2}\lesssim |E|^\frac12|F|^\frac12$.
Assuming without loss of generality that $|E| \leq |F|$,  we have
\begin{equation} \label{doublingnew} {||B_{\theta}(\chi_E,\chi_F)||}_{L^2({\Bbb R}^2)}^2
=\int \int \int \chi_E(x-y) \chi_F(x-\Theta y)\chi_E(x-y') \chi_F(x-\Theta y')\, d\sigma(y) d\sigma(y') dx \end{equation}
$$=\int_{|\alpha-\alpha'|<\frac{\theta}{2}} \int \int \chi_E(x-y) \chi_F(x-\Theta y)\chi_E(x-y') \chi_F(x-\Theta y')\, d\sigma(y) d\sigma(y') dx $$
$$\qquad+ \int_{|\alpha-\alpha'| \ge \frac{\theta}{2}} \int \int \chi_E(x-y) \chi_F(x-\Theta y)\chi_E(x-y') \chi_F(x-\Theta y')\, d\sigma(y) d\sigma(y')dx $$
$$ \leq \int_{|\alpha-\alpha'|<\frac{\theta}{2}}  \int \int \chi_E(x-y) \chi_F(x-\Theta y')\, d\sigma(y) d\sigma(y') dx\qquad\qquad\qquad\qquad\qquad$$
$$+\int_{|\alpha-\alpha'| \ge \frac{\theta}{2}} \int \int  \chi_E(x-y) \chi_E(x- y')\, d\sigma(y) d\sigma(y') dx,$$ where $y=(\cos(\alpha), \sin(\alpha))$, $y'=(\cos(\alpha'), \sin(\alpha'))$. 

\vskip.125in 

To make a change of variables for the first integral in the last expression, we  consider 
\begin{equation} \label{cov} u_1=x_1-\cos(\alpha), \ u_2=x_2-\sin(\alpha), \ v_1=x_1-\cos(\alpha'+\theta), v_2=x_2-\sin(\alpha'+\theta);
 \end{equation}
for the second integral, we make the change of variables 
\begin{equation}
u_1=x_1-\cos(\alpha), \ u_2=x_2-\sin(\alpha), \ v_1=x_1-\cos(\alpha'), v_2=x_2-\sin(\alpha').
\end{equation}
The Jacobian for the first  is 
$$ \sin(\alpha-\alpha'-\theta),$$ 
while the Jacobian of the  second is 
$$ \sin(\alpha-\alpha').$$  
Note that both of these quantities are bounded away from $0$ because of the constraints  on the angle between $y$ and $y'$. Note that this argument fails when $\theta=\pi$ since if $|\alpha-\alpha'| \geq \frac{\theta}{2}$, the Jacobian goes to $0$ regardless of which terms we keep.
\medskip

As long as $0<\theta<\pi$, though, we have that the Jacobian in both cases is bounded from below by $\frac{1}{2}sin(\frac{\theta}{2})$. It follows that 
$$ {||B_{\theta}(\chi_E,\chi_F)||}_{L^2({\Bbb R}^2)} \leq C_{\theta}(|E|^2+|E||F|)^{\frac{1}{2}} \leq 2C {|E|}^{\frac{1}{2}} {|F|}^{\frac{1}{2}}.$$ for some constant $C$ depending only on $\theta$. Moreover, it is not difficult to see from the argument above that 
$$ C_{\theta} \leq C' \frac{1}{\min \{\theta, \pi-\theta \}},$$ where $C'$ is a uniform constant independent of $\theta$, 
establishing boundedness of $B_\theta$ at the vertex (vii)  for nondegenerate $\theta$, i.e., $\theta\in (0,\pi)$. 

\vskip.25in

\section{The $L^{\frac{3}{2}} \times L^{\frac{3}{2}} \to L^1$ estimate for general operators} 
\label{sec L1 general}

\vskip.125in

We now show that, if (i) $(\phi_1(x,y),\phi_2(x,z),\phi_3(y,z))$ satisfy  (\ref{smooth Z}), 
and
(ii) $\phi_3(y,z)$ satisfies the Phong-Stein condition (\ref{cond PS}) for all $(y,z)\in Z_3$, then
$B:L^\frac32(\R^2)\times L^\frac32(\R^2) \to L^1(\R^2)$.
As for the model operators, since $K\ge 0$,  one can assume that $f,g\ge 0$, and write
$$||B(f,g)||_{L^1}=\int B(f,g)(x) dx=\int f(y)\, \left[\int \left[ \int K(x,y,z)\, dx\right] g(z)\, dz\right]\, dy.$$
The expression inside the outermost of the square brackets is a linear operator, $Tg(y)$. 
It hence suffices to show that $T:L^\frac32(\R^2)\to L^3(\R^2)$;  by the Phong-Stein condition (\ref{cond PS}) for $\phi_3$, 
this will follow if we show that  $T$ is a generalized Radon transform of associated to $Z_3$. 
The kernel of $T$ is
$$L(y,z)=\int K(x,y,z)\, dx = \pi_*(K)(y,z),$$
where $\pi_*$ denotes pushforward of distributions under the projection $\pi(x,y,z)=(y,z)$. 
The operator $\pi_*:\mathcal E'(\R^6)\to \mathcal E'(\R^4)$  is itself an FIO, $\pi_*\in I^{-\frac12}(C_\pi)$, associated to the canonical relation
$$C_\pi=\left\{(y,\eta,z,\zeta; x,y,z,0,\eta,\zeta): (x,y,z)\in\R^6,\, (\eta,\zeta)\in\R^4\setminus 0\right\};$$
see Guillemin and Sternberg \cite{GS}.
$C_\pi$ is a nondegenerate canonical relation in $T^*\R^4\times T^*\R^{10}$,  and thus its application to $N^*Z$ is covered by the transverse intersection calculus. A direct calculation shows that $C_\pi\circ N^*Z=N^*Z_3$. 
 and thus $L\in I^{-\frac12}(N^*Z_3)$. Hence,  $T$  is a  linear generalized Radon transform on $\R^2$ satisfying the Phong-Stein condition, 
 and has the same mapping properties as  any such operator; in particular,  $T:L^\frac32(\R^2)\to L^3(\R^2)$. 
\medskip

\section{Proof of the $L^{2,1}\times L^{2,1}\to L^2$  estimate for general operators} 
\label{sec L2 general}

\vskip.125in 

Now, to prove the restricted strong type $L^{2,1}({\Bbb R}^2) \times L^{2,1}({\Bbb R}^2) \to L^2({\Bbb R}^2)$  result for $B$,  consider as in Sec. \ref{sec L2 model} the $L^2$ norm squared of the operator applied to indicator functions of sets $E, \, F$: 
\begin{equation} \label{squaredthing} 
\int \int \int \int \int \chi_E(y)\chi_E(y')\chi_F(z)\chi_F(z') K(x,y,z)K(x,y',z')\, dydy'dzdz'dx. 
\end{equation}

Modifying somewhat  the argument for the model operators in Sec. \ref{sec L2 model},
one can show that if any one of the following four  bounded properties holds,  one  obtains $||B(\chi_E,\chi_F)||_{L^2}\lesssim |E|^\frac12|F|^\frac12$: 

\be\label{Kyz'}
K_{yz'}:=\int \int \int K(x,y,z)K(x,y',z')\, dy'dzdx\in L^\infty(\R^6_{x,y,z'}),
\ee 

\be\label{Kzy'}
K_{zy'}:=\int \int \int K(x,y,z)K(x,y',z')\, dydz'dx\in L^\infty(\R^6_{x,z,y'}),
\ee

\be\label{Kzz'}
K_{zz'}:=\int \int \int K(x,y,z)K(x,y',z')\, dydy'dx\in L^\infty(\R^6_{x,z,z'}),
\ee
and

\be\label{Kyy'}
K_{yy'}:=\int \int \int K(x,y,z)K(x,y',z')\, dz dz' dx\in L^\infty(\R^6_{x,y,y'}).
\ee

More precisely, if (\ref{Kyz'}) holds, we eliminate $z,\, y'$  in (\ref{squaredthing}) by noting  that $\chi_E(y')\chi_F(z)\le 1$  and obtain 
an  upper bound $C|E||F|$  for (\ref{squaredthing}). 
If (\ref{Kzy'}) holds, we proceed in a similar way, using $\chi_E(y)\chi_F(z')\le 1$. 
If (\ref{Kzz'}) holds, we employ $\chi_E(y)\chi_E(y')\le 1$ and bound the whole expression in (\ref{squaredthing}) by $C{|F|}^2$,
which, if $|F|\le|E|$, is bounded by $C|E||F|$. 
On the other hand, if $|E|\le|F|$, we may use the boundedness of  the expression in (\ref{Kyy'}), together with $\chi_F(z)\chi_F(z')\le 1$, 
to bound (\ref{squaredthing}) by $C{|E|}^2$, which is $\leq C|E||F|$ . 
Thus, regardless of whether $|E|\le |F|$ or $|F|\le |E|$, we obtain
$||B(\chi_E,\chi_F)||_{L^2}\lesssim |E|^\frac12|F|^\frac12$.  
This argument holds more generally if 
any  $(x,y,z,y',z')$ has a neighborhood in $\R^{10}$ on which  one of 
\be\label{K logic}
\hbox{(\ref{Kyz'}) or (\ref{Kzy'}) or [ (\ref{Kzz'}) and (\ref{Kyy'}) ]  holds,}
\ee
with the conjunction in the last term to cover both of the cases $|E|\le |F|$ and $|F|\le |E|$.
Taking a subordinate  partition of unity of $\R^{10}$, the domain of integration  in (\ref{squaredthing}), we can then apply the above arguments to still obtain
$||B(\chi_E,\chi_F)||_{L^2}\lesssim |E|^\frac12|F|^\frac12$.

\vskip.125in

In the framework of the product co-normal kernels above, we can formulate first order conditions  on the $\phi_j$, and hence on $Z$ and $K$, 
which imply that one of   (\ref{Kyz'})-(\ref{Kyy'}) holds locally. 
Let
$$\hKK(x,y,z,y',z')=K(x,y,z)\cdot K(x,y',z')\in\mathcal D'(\R^{10})$$
and
$$\hZZ=\left\{(x,y,z,y',z'):(x,y,z)\in Z,\, (x,y',z')\in Z\right\}\subset\R^{10}.$$
Analogous to (\ref{smooth Z}), we assume
\be\label{smooth ZZ}
 rank\left[
\begin{matrix}
d_x\phi_1(x,y) & d_y\phi_1& 0& 0 & 0 \cr
d_x\phi_2(x,z) & 0 & d_z\phi_2 & 0 & 0\cr
0 & d_y\phi_3(y,z) & d_z \phi_3 & 0 & 0 \cr
d_x\phi_1(x,y') & 0 & 0 & d_{y'}\phi_1& 0 \cr
d_x\phi_2(x,z')& 0 & 0 & 0 & d_{z'}\phi_2\cr
0 & 0 & 0 &d_{y'}\phi_3(y',z') & d_{z'} \phi_3
\end{matrix}
\right]=6,
\ee
for all $(x,y,z,y',z')\in\hZZ$. Then $\hZZ$ is smooth and 4-dimensional, and $\hKK$ is a smooth density on it. 
The kernel $K_{yz'}$ in (\ref{Kyz'}) above is just $\pi^{yz'}_*(\hKK)(y,z')\in\mathcal D'(\R^4)$, the pushforward by  $\pi^{yz'}(x,y,z,y',z')=(y,z')$. This will be a smooth density on $\R^4$, 
hence with a smooth (thus locally bounded) Radon-Nikodym derivative, if $\hZZ$ is a graph over the $y,z'$ variables, 
and by the implicit function theorem this holds when the $6\times 4$ submatrix consisting 
of the $y$ and $z'$ columns of the $6\times 10$ matrix in (\ref{smooth ZZ}) has rank 4. Deleting the two rows of zeros, 
this holds at a point $(x_0,y_0,z_0,y_0',z_0')\in\hZZ$  iff

\be\label{cond yz'}
 det \left[
\begin{matrix}
d_y\phi_1(x_0,y_0) & 0 \cr
d_y\phi_3(y_0,z_0) & 0 \cr
0 & d_{z'}\phi_3(x_0,z_0')  \cr
0 & d_{z'}\phi_3(y_0',z_0') 
\end{matrix}
\right]\ne 0.
\ee
Similarly, considering the pushfoward maps $\pi_*^{zy'},\, \pi_*^{zz'}$ and $\pi_*^{yy'}$, 
we see that the local  boundedness of (\ref{Kzy'}), (\ref{Kzz'}) and (\ref{Kyy'}) follow from

\be\label{cond zy'}
 det \left[
\begin{matrix}
d_z\phi_2(x_0,z_0) & 0 \cr
d_z\phi_3(y_0,z_0) & 0 \cr
0 & d_{y'}\phi_1(x_0,y_0')  \cr
0 & d_{y'}\phi_3(y_0',z_0') 
\end{matrix}
\right]\ne 0,
\ee

\be\label{cond zz'}
 det \left[
\begin{matrix}
d_z\phi_2(x_0,z_0) & 0 \cr
d_z\phi_3(y_0,z_0) & 0 \cr
0 & d_{z'}\phi_2(x_0,y_0')  \cr
0 & d_{z'}\phi_3(y_0',z_0') 
\end{matrix}
\right]\ne 0,
\ee
and

\be\label{cond yy'}
 det \left[
\begin{matrix}
d_y\phi_1(x_0,y_0) & 0 \cr
d_y\phi_3(y_0,z_0) & 0 \cr
0 & d_{y'}\phi_1(x_0,y_0')  \cr
0 & d_{y'}\phi_3(y_0',z_0') 
\end{matrix}
\right]\ne 0,\quad\hbox{resp.}
\ee

\bigskip

Conditions (\ref{cond yz'})-(\ref{cond yy'}) are of course open conditions, so if any one holds at a point $(x_0,y_0,z_0,y_0',z_0')\in\hZZ$, 
then it holds on a neighborhood. Thus, if at every $(x_0,y_0,z_0,y_0',z_0')\in\hZZ$,

\be\label{cond general}
\hbox{at least one of (\ref{cond yz'}) or (\ref{cond zy'}) or [ (\ref{cond zz'}) and (\ref{cond yy'}) ] holds,}
\ee
then (\ref{K logic}) holds. Combined with the discussion in Sec. \ref{sec L1 general}, valid if $\phi_3$ 
satisfies the rotational curvature condition, this finishes the proof of Thm. \ref{thm general}.
\medskip

\vskip.125in

\end{document}